\documentclass[12pt]{article}
\usepackage{amsmath}
\usepackage{amssymb}
\usepackage{cite}

\newcommand{\eh}{\hfill}\newlength{\sperr}

\newenvironment{proof}{{\settowidth{\sperr}{\bf\rm
Proof}%
\par\addvspace{0.3cm}\noindent\parbox[t]{1.3\sperr}
{\bf\rm P\eh r\eh o\eh o\eh f\eh }%
}}{\nopagebreak\mbox{}\hfill
$\Box$\par\addvspace{0.3cm}}

\def\nn{\nonumber}
\def\a{\alpha}
\def\b{\beta}
\def\g{\gamma}
\def\G{\Gamma}

\def\vk{\varkappa}

\def\om{\omega}

\def\t{\theta}

\def\vp{\varphi}

\def\ve{\varepsilon}
\def\wh{\widehat}
\def\wt{\widetilde}
\def\ov{\overline}

\def\BC{{\mathbb C}}
\def\BR{{\mathbb R}}

\def\clp{{\mathcal P}}

\def\cln{{\mathcal N}}
\def\clu{{\mathcal U}}

\def\im{{\rm Im\ }}

\newcommand{\E}{\mathrm{e}}
\newcommand{\I}{\mathrm{i}}

\newtheorem{Pa}{Paper}[section]
\newtheorem{Tm}[Pa]{{\bf Theorem}}
\newtheorem{La}[Pa]{{\bf Lemma}}
\newtheorem{Cy}[Pa]{{\bf Corollary}}
\newtheorem{Rk}[Pa]{{\bf Remark}}

\newtheorem{Dn}[Pa]{{\bf Definition}}
\newtheorem{Pn}[Pa]{{\bf Proposition}}

\title{Recovery of Dirac system from the  rectangular 
Weyl matrix function}

\author{B. Fritzsche, B. Kirstein, I.Ya. Roitberg, A.L. Sakhnovich}

\date{}
\begin{document}
\maketitle

\begin{abstract} Weyl theory  for
Dirac systems with rectangular matrix potentials
is non-classical. The corresponding Weyl functions are
rectangular matrix functions. Furthermore, they are non-expansive in the 
upper semi-plane.
Inverse problems are treated for such Weyl
functions, and some results are new even for the  square
Weyl functions. High energy asymptotics of Weyl functions
and Borg-Marchenko type uniqueness results are derived too.
\end{abstract}

{MSC(2010):}  34B20, 34L40.

Keywords:  {\it Weyl function, Weyl theory, Dirac system, rectangular matrix potential,  inverse problem.}

\section{Introduction} \label{intro}
\setcounter{equation}{0}
 
The self-adjoint Dirac-type (also called Dirac, ZS or AKNS)
system
\begin{equation}       \label{1.1}
\frac{d}{dx}y(x, z )=\I (z j+jV(x))y(x,
z ) \quad
(x \geq 0),
\end{equation}
where
\begin{equation}   \label{1.2}
j = \left[
\begin{array}{cc}
I_{m_1} & 0 \\ 0 & -I_{m_2}
\end{array}
\right], \hspace{1em} V= \left[\begin{array}{cc}
0&v\\v^{*}&0\end{array}\right],
 \end{equation}
$I_{m_k}$ is the $m_k \times
m_k$ identity
matrix and $v(x)$ is an $m_1 \times m_2$ matrix function,
is a classical matrix differential equation. We always assume
that  $v$ is measurable and, moreover, locally summable, that is, summable
on all the finite intervals.

Dirac-type systems   are very well-known in
mathematics and
applications, especially in mathematical physics (see, e.g., books
\cite{AD, CS, K, LS, SaL3},
recent publications \cite{AGKLS1, AGKLS2, AD0,  CG1, CG2, FKS1, FKS2, SaL4},
and numerous references therein). System \eqref{1.1} is used, in particular,
in the study of transmission lines and acoustic problems \cite{YL}.
The most interesting applications are, however, caused by the fact that system 
(\ref{1.1}) is an auxiliary linear system for many  important nonlinear integrable
wave equations and as such it was studied, for instance, in
\cite{AKNS, AbPT, AS, FT, GSS, Kaup, KS2, SaA7, ZS}. (The
name ZS-AKNS system is connected with these applications.) Nonlinear
Schr\"odinger equations, modified Korteweg de Vries equations, and
second harmonic generation model, which describe various wave processes (including, e.g.,
nonlinear wave processes in water and in waveguides, in nonlinear optics, photoconductivity, and on silicon
surfaces),  are only some of the well-known examples.
The evolution of Weyl function for these systems is described in terms of
M\"obius transformations (see, e.g., \cite{SaA1, SaA7, ALS11, SaLev1, SaLev2, SaL3}
and references therein), which is one of the fruitful approaches
to study interesting initial-boundary value problems for integrable nonlinear
equations.

The Weyl and spectral theory of self-adjoint Dirac systems, where $m_1=m_2$,
was treated, for instance, in \cite{AD0, Kre1, CG1, LS, SaA3,  SaL3}
(see also various references therein). The "non-classical" Weyl theory for the
equally important  case $m_1\not= m_2$, which appears in the study of  coupled, multicomponent, and matrix nonlinear equations,
is the subject of this paper. 

The $m_1\times m_2$ matrix function $v(x)$ from \eqref{1.2} is called the potential of system
\eqref{1.1}. We put $m_1+m_2=:m$. The fundamental solution of system
\eqref{1.1} is denoted by $u(x,z)$, and this solution is normalized by the condition
\begin{align}&      \label{1.3}
u(0,z)=I_m.
\end{align}

The next Section \ref{Snode} is dedicated to representation of the fundamental solution.
In Section \ref{HE} we follow \cite{FKRSp1} to introduce the Weyl function,
and study the high energy asymptotics of this Weyl function. The solution of the inverse
problem to recover the potential (and system) from the Weyl function is given in Section
\ref{IP}. Borg-Marchenko-type uniqueness results are contained in that section as well.
Finally, Section \ref{Sc} is dedicated to conditions for an analytic  matrix function
to be the Weyl function of some Dirac-type system. The results from this section
were not published before even for the case that $m_1=m_2$, though some of them appeared for that case
in the thesis \cite{SaAth}.

As usual, $\BC$ stands for the complex plain,
$\BC_+$ for the open upper
semi-plane, and $\im$ stands for image. If $a\in \BC$, then $\ov a$ is its complex conjugate. 
The notation $L^2_{m_2\times m_1}(0,\, \infty)$ will be used to denote the space of $m_2 \times m_1$
matrix-functions with entries belonging to $L^2(0,\, \infty)$. 
We use $\ov{\bf H}$ to denote the closure of the space ${\bf H}$, $I$ to denote
the identity operator, and
$B({\bf H_1}, {\bf H_2})$ to denote the
class of bounded operators acting from ${\bf H_1}$ to ${\bf H_2}$. We write 
$B({\bf H_1})$ if ${\bf H_1}={\bf H_2}$, and $B^1[0,\,l]$ will be used to denote
a class of functions (or matrix functions),  whose derivatives are bounded on 
$[0,\,l]$. 
An $m_2 \times m_1$ matrix $\a$ is said to be non-expansive, if $\a^*\a \leq I_{m_1}$ (or, equivalently, if $\a\a^* \leq I_{m_2}$).

\section{Representation of fundamental solution}\label{Snode}
\setcounter{equation}{0}
The results of this section can be formulated for Dirac system on a fixed final interval
$[0,\, l]$. We assume that $v$ is bounded on this interval and put
\begin{align}&      \label{3.1}
\b(x)=\begin{bmatrix}
I_{m_1} & 0
\end{bmatrix}u(x,0), \quad \g(x)=\begin{bmatrix}
0 &I_{m_2}
\end{bmatrix}u(x,0).
\end{align}
It follows from $\sup_{x < l}\|v\|<\infty$ and from \eqref{1.1} that
\begin{align}&      \label{3.2}
\sup_{x < l}\|\g^{\prime}(x)\|<\infty, \quad \g^{\prime}:=\frac{d}{dx}\g .
\end{align}
Moreover, by \eqref{1.3} and \eqref{2.4} we get $u(x,0)^*ju(x,0)=j=u(x,0)ju(x,0)^*$. Therefore, 
\eqref{3.1} implies
\begin{align}&      \label{3.3}
\b j \b^*\equiv I_{m_1}, \quad \g j \g^*\equiv -I_{m_2}, \quad \b j \g^*\equiv 0. 
\end{align}

Next, we need the following similarity result for the Volterra operator
\begin{align}&      \label{3.4}
K=\int_0^x  \, F(x)G(t) \cdot \, dt, \quad K \in B\big(L^2_{m_1}(0, \, l)\big),
\end{align}
where $F(x)$ is an $m_1 \times m$ matrix function, $G(t)$ is an $m \times m_1$ matrix function,
and
\begin{equation} \label{3.5}
\I F(x)G(x) \equiv I_{m_2}.
\end{equation}
\begin{Pn}\label{PnSim}\cite{SaL0} Let $F$ and $G$ be boundedly differentiable
and let (\ref{3.5}) hold. Then we have
\begin{equation} \label{3.6}
K= EAE^{-1}, \quad A:=-\I \int_0^x \, \cdot \, dt,  \quad A, E, E^{-1} \in B\big(L^2_{m_2}(0, \, l)\big),
\end{equation}
where $K$ is given by (\ref{3.4}) and $E$ is a triangular operator of the form
\begin{equation} \label{3.7}
(E f)(x)=\rho(x) f(x)+\int_0^x E(x,t)f(t)dt, \quad \frac{d}{d x}\rho =\I F^{\prime}G\rho, \quad
\det \rho(0)\not=0.
\end{equation}
Moreover, the operators $E^{\pm 1}$ map functions with bounded derivatives
into functions with bounded derivatives.
\end{Pn}
The proposition above is a particular case of Theorem 1 \cite{SaL0}
(see also a later paper \cite{AGKLS2} on the case of continuously
differentiable $F$ and $G$). 

Set
\begin{align}&      \label{3.8}
F=\g, \quad G=\I j\g^*, \quad \g j \g^*\equiv -I_{m_2}, \quad \g \in B^1[0,\,l],
\end{align}
where $B^1$ stands for the class of boundedly differentiable matrix functions. 
Clearly $F$ and $G$ in \eqref{3.8} satisfy conditions of Proposition \ref{PnSim}. 
Separate $\g(x)$ into two blocks 
$\g=\begin{bmatrix}
\g_1 & \g_2
\end{bmatrix}$, where $\g_1, \, \g_2$ are $m_2\times m_1$ and $m_2\times m_2$,
respectively, matrix functions. 
Without loss of generality one can choose $E$ so that 
\begin{align}&      \label{3.12}
E^{-1}\g_2\equiv I_{m_2},
\end{align}
where $E^{-1}$ is applied to $\g_2$ columnwise.
\begin{Pn}\label{PnSimN}Let 
$K$ be given by \eqref{3.4}, where $F$ and $G$ satisfy \eqref{3.8}, and let $\wt E$
be a similarity operator from Proposition \ref{PnSim}.
Introduce $E_0 \in B\big(L^2_{m_2}(0, \, l)\big)$ by the equalities
\begin{align}&      \label{3.13}
\big(E_0f\big)(x)=\rho(0)^{-1}\g_2(0)f (x) + \int_0^x E_0(x-t) f(t) dt, \quad 
E_0(x):=\big(\wt E^{-1}\g_2\big)^{\prime}(x).
\end{align}
Then, the operator $E:=\wt E E_0$ is another similarity operator from Proposition \ref{PnSim},
which satisfies the additional condition   \eqref{3.12}.
\end{Pn}
\begin{proof}. 
The proof of the proposition is similar to the case $m_1=m_2$ (see, e.g.,
\cite[pp. 103, 104]{SaL0}). Indeed, the next identity can be easily shown directly 
(and follows also from the fact that $E_0$ is a lower triangular convolution operator):
\begin{align}&      \label{p1}
A E_0=E_0A, \quad ( A=-\I \int_0^x \, \cdot \, dt).
\end{align}
Furthermore, because of the third relation in \eqref{3.8} we have
$\det \g_2(0)\not=0$, and so $E_0$ is invertible. Hence, equalities \eqref{p1}
and $K= \wt EA \wt E^{-1}$ imply \eqref{3.6}:
\begin{align}&      \label{p2}
K= EA  E^{-1}, \quad E:=\wt E E_0.
\end{align}
Formula \eqref{3.13} leads us also to the equality
\begin{align}     \nn
\big(E_0 I_{m_2}\big)(x)&=\rho(0)^{-1}\g_2(0)+\int_0^xE_0(x-t)dt=
\rho(0)^{-1}\g_2(0)+\int_0^xE_0(t)dt
\\  \label{p3}&
=\rho(0)^{-1}\g_2(0)+
\big(\wt E^{-1} \g_2 \big)(x)-\big(\wt E^{-1} \g_2 \big)(0).
\end{align}
Recalling that $\wt E^{-1} \g_2 \in B^1[0,\,l]$ we obtain a representation:
\begin{align}&      \label{p4}
\big(\wt E^{-1} \g_2 \big)(x)=\rho(0)^{-1}\g_2(0)+\I\Big(A
\big(\wt E^{-1} \g_2 \big)^{\prime}\Big)(x).
\end{align}
Using \eqref{p4} we rewrite \eqref{p3} as
\begin{align}     \label{p5}&
\big(E_0 I_{m_2}\big)(x)=
\big(\wt E^{-1} \g_2 \big)(x),
\end{align}
and \eqref{3.12} follows. It remains to show
that $E_0^{\pm 1}$ maps $B^1[0,\,l]$ into $B^1[0,\,l]$. 
First note that the integral operators $E_0, \,E_0^{-1}$ have  bounded kernels
and map bounded functions into bounded. In particular,  for  $E_0^{-1}$ it follows
from a series representation of the operator of the form
$\Big(I+\int_0^x k(x,t) \, \cdot \, dt\Big)^{-1}$
(where the kernel $k$ is bounded). Now, taking into account
that (similar to \eqref{p4}) any $f \in B^1[0,\,l]$ admits representation
$f=f(0)+\I Af^{\prime}$, we see that formulas \eqref{p1}
and \eqref{p5} yield the fact that $E_0$ maps $B^1[0,\,l]$ into $B^1[0,\,l]$.
The identity $A E_0^{-1}=E_0^{-1}A$ is immediate from \eqref{p1}.
Thus,  to prove that $E_0^{-1}$ maps $B^1[0,\,l]$ into $B^1[0,\,l]$ we need
only to show that $E_0^{-1}I_{m_2}\in B^1[0,\,l]$, which can be derived
from \eqref{p4} and \eqref{p5}:
\begin{align}     \nn&
E_0^{-1} I_{m_2}=
\Big(I_{m_2}-\I E_0^{-1}A\big(\wt E^{-1}\g_2\big)^{\prime}\Big)\g_2(0)^{-1}\rho(0)
 \in B^1[0,\,l].
\end{align}
\end{proof}
\begin{Rk}\label{kern} The kernels of the operators $\wt E^{\pm 1}$, which are
constructed in \cite{SaL0}, as well as the kernels of the operators $E_0^{\pm 1}$
from the proof of Proposition \ref{PnSimN} are bounded.
Therefore, without loss of generality we always assume further that the kernels of $E^{\pm 1}$ are bounded.
\end{Rk}
The next lemma easily follows from Proposition \ref{PnSimN}
and will be used to construct fundamental solution.
\begin{La}\label{LaNode} Let $\g$ be an $m_2 \times m$ matrix function, which satisfies the last two relations in \eqref{3.8}, and set
\begin{align}&      \label{p6}
S:=E^{-1}\big(E^*\big)^{-1}, \quad \Pi:= \begin{bmatrix}
\Phi_1 & \Phi_2
\end{bmatrix}, \quad 
\Phi_k \in B\big(\BC^{m_k}, \, L^2_{m_2}(0, \, l)\big);
\\
&      \label{3.19}
\big(\Phi_1 f\big)(x)=\Phi_1(x)f, \quad
\Phi_1(x):=\big(E^{-1}\g_1\big)(x);  \quad  \Phi_2 f=I_{m_2}f\equiv f;
\end{align}
where $E$ is constructed (for the given $\g$) in Proposition \ref{PnSimN}.
Then $A$, $S$, and $\Pi$ form an $S$-node, that is (see \cite{SaL1, SaL2, SaL3}),
the operator identity
\begin{align}&      \label{p9}
AS-SA^*=\I \Pi j \Pi^*
\end{align}
holds. Furthermore,
we have
\begin{align}&      \label{p7}
\ov{\sum_{i=0}^{\infty}\im \Big(\big(A^*\big)^i S^{-1}\Pi\Big)}=L^2_{m_2}(0, \, l).
\end{align}
\end{La}
\begin{proof}.  Because of \eqref{3.4},  \eqref{3.6} and \eqref{3.8}, we get
\begin{align}&      \label{p8}
EAE^{-1}-\big(E^{-1}\big)^*A^*E^*=K-K^*=\I\g(x)j\int_0^l\g(t)^*\,\cdot \,dt.
\end{align}
Formulas \eqref{3.12},  \eqref{p6} and \eqref{3.19} lead us to the equality
\begin{align}&      \label{p10}
\Pi f=\big(E^{-1}\g\big)(x)f.
\end{align}
Now, the operator identity \eqref{p9} follows from \eqref{p8},  \eqref{p10}, and the first equality in \eqref{p6}.

To prove \eqref{p7} we will show that
\begin{align}&      \label{p11}
\sum_{i=0}^{N}\im \Big(\big(A^*\big)^i S^{-1}\Pi\Big)
\supseteq \sum_{i=0}^{N}\im \big(S^{-1}A^i\Pi\big)
=S^{-1}\sum_{i=0}^{N}\im \big(A^i\Pi\big).
\end{align}
For that purpose we rewrite \eqref{p9} as $S^{-1}A=A^*S^{-1}+\I
S^{-1}\Pi j \Pi^* S^{-1}$. Hence, for $N_1, N_2 \geq 0$ we obtain
\begin{align}&     \label{p12}
\im \Big(\big(A^*\big)^{N_1+1} S^{-1}A^{N_2}\Pi\Big)
\\ \nn &
+
\sum_{i=0}^{N_1+N_2}\im \Big(\big(A^*\big)^i S^{-1}\Pi\Big)
\supseteq \im \Big(\big(A^*\big)^{N_1} S^{-1}A^{N_2+1}\Pi\Big).
\end{align}
Using \eqref{p12}, we derive \eqref{p11} by induction.
In view of \eqref{p11}, it suffices to show that
\begin{align}&      \label{p13}
\ov{\sum_{i=0}^{\infty}\im \big(A^i \Pi\big)}=L^2_{m_2}(0, \, l),
\end{align}
which, in its turn, follows from \eqref{3.19}.
\end{proof}
\begin{Rk} \label{TrMF}
Given an $S$-node \eqref{p9}, we introduce a transfer matrix function
in Lev Sakhnovich form $($see \cite{SaL1, SaL2, SaL3}$):$
\begin{align}&      \label{p13a}
w_A(r,z):=I_m+\I zj\Pi^*S_r^{-1}(I-zA_r)^{-1}P_r\Pi, \quad 0<r\leq l,
\end{align}
where $I$ is the identity operator; $A_r, \, S_r \, \in \, B\big(L^2_{m_2}(0, \, r)\big)$, 
\begin{align}&      \label{p13b}
A_r:=P_rAP_r^*, \quad S_r:=P_rSP_r^*;
\end{align}
 $A$ is given by \eqref{3.6}, the operators $S$ and $\Pi$ are given by
 \eqref{p6} and \eqref{3.19}, and the
operator $P_r$ is an orthoprojector from $L^2_{m_2}(0, \, l)$ on $L^2_{m_2}(0, \, r)$ such that
\begin{align}&      \label{p13c}
\big(P_rf\big)(x)=f(x) \quad (0<x<r), \quad f \in L^2_{m_2}(0, \, l).
\end{align}
Since $P_rA=P_rAP_r^*P_r$, it follows from \eqref{p9} that the operators $A_r$, $S_r$, and $P_r\Pi$ form an 
$S$-node too, that is, the operator identities
\begin{align}&      \label{p9'}
A_rS_r-S_rA_r^*=\I P_r \Pi j \Pi^*P_r^*
\end{align}
hold.
\end{Rk}
Now, in a way similar to \cite{SaA0, SaA1} the fundamental solution $w$
of the system
\begin{align}&      \label{3.14}
\frac{d}{dx}w(x,z)=\I zj\g(x)^*\g(x)w(x,z), \quad w(0,z)=I_m
\end{align}
is constructed.
\begin{Tm}\label{FundSol}
Let $\g$ be an $m_2 \times m$ matrix function, which satisfies the last two relations in \eqref{3.8}.
 Then, the fundamental solution $w$ given by \eqref{3.14} 
 admits representation
 \begin{align}&      \label{3.16}
w(r,z)=w_A(r,z),
\end{align}
where $w_A(r,z)$ is defined in Remark \ref{TrMF}.
\end{Tm}
\begin{proof}. The statement of the theorem follows from Continual factorization theorem
(see \cite[p. 40]{SaL3}). More precisely, our statement follows from a corollary of the
Continual factorization theorem, namely, from Theorem 1.2 \cite[p. 42]{SaL3}.
Using Lemma \ref{LaNode} we easily check that the conditions 
of Theorem 1.2 \cite[p. 42]{SaL3} are fulfilled. Therefore, if $\Pi^*S_r^{-1}P_r\Pi$
is boundedly differentiable, we have
\begin{align}&      \label{p14}
\frac{d}{dr}w_A(r,z)=\I zjH(r)w_A(r,z), \quad \lim_{r\to +0}w_A(r,z)=I_m, \\
&      \label{p14'}
 H(r):=\frac{d}{dr}\big(\Pi^*S_r^{-1}P_r\Pi\big),
\end{align} 
where $w_A$ is given by \eqref{p13a}. Since $E^{\pm 1}$ are lower triangular operators, we see that
\begin{align}&      \label{p14!}
P_rEP_r^*P_r=P_rE, \qquad  \big(E^{-1}\big)^*P_r^*=P_r^*P_r\big(E^{-1}\big)^*P_r^*.
\end{align} 
Hence, formulas \eqref{p6} and
\eqref{p13b} lead us to
\begin{align}&      \label{p14!!}
 S_r^{-1}=E_r^*E_r, \qquad E_r:=P_rEP_r^*.
\end{align} 
Therefore, taking into account  \eqref{p10}, we rewrite  \eqref{p14'} as
\begin{align}&      \label{p15}
 H(r)=\g(r)^*\g(r).
\end{align} 
Formulas \eqref{3.14},
\eqref{p14} and \eqref{p15} imply \eqref{3.16}.
\end{proof}
Now, consider again the case of Dirac system.  Because of  \eqref{3.1} and \eqref{3.3} we obtain
\begin{align}&      \label{p16}
u(x,0)j\g(x)^*\g(x)u(x,0)^{-1}=-\begin{bmatrix}
0 & 0\\0 & I_{m_2}
\end{bmatrix}.
\end{align}
Hence, direct calculation shows that the following corollary of Theorem \ref{FundSol}
is true.
\begin{Cy}\label{FSD} Let $u(x,z)$ be the fundamental solution of a Dirac system with a bounded potential
$v$ and let $\g$ be given by \eqref{3.1}. Then $u(x,z)$ admits representation
\begin{align}&      \label{p17}
u(x,z)=\E^{ixz}u(x,0)w(x,2z),
\end{align}
 where $w$ has the form \eqref{3.16} and the $S$-node generating the transfer matrix function $w_A$    
is recovered from $\g$ in Lemma \ref{LaNode}.
\end{Cy}
\begin{Rk}\label{beta'}
For the case that $\g$ is given by \eqref{3.1}, it follows from   \eqref{1.1} and \eqref{3.3} that
\begin{align}&      \label{3.9}
\g^{\prime}(x)j\g(x)^*=-\I \begin{bmatrix}
v(x)^* & 0
\end{bmatrix}u(x,0)j\g(x)^*=-\I v(x)^*\b(x)j\g(x)^*\equiv 0.
\end{align}
Thus, from \eqref{3.8} and \eqref{3.9} we see that $\rho$ in \eqref{3.7} is a constant matrix.
Therefore, since $\g_2(0)=I_{m_2}$, equality \eqref{p10} implies that  $\rho(x)\equiv I_{m_2}$, and 
formula \eqref{3.7} can be rewritten in the form
\begin{equation} \label{3.10}
(E f)(x)=f(x)+\int_0^x E(x,t)f(t)dt, \quad E \in B\big(L^2_{m_2}(0, \, l)\big).
\end{equation}
\end{Rk}
Recalling that $\det \g_2(x)\not=0$ we can rewrite \eqref{3.9} as
\begin{align}&      \label{p18}
\g_2^{\prime}=\g_1^{\prime}(\g_2^{-1}\g_1)^*.
\end{align}
Using \eqref{p18}, we recover $\g_2$ and $\g_1=\g_2(\g_2^{-1}\g_1)$
from $\g_2^{-1}\g_1$. Indeed, we have
\begin{align}&      \nn
(\g_2^{-1}\g_1)^{\prime}=-\g_2^{-1}\g_2^{\prime}\g_2^{-1}\g_1+\g_2^{-1}
\g_1^{\prime}, \,\,{\mathrm{i.e.}},\,\, \g_1^{\prime}=\g_2^{\prime}(\g_2^{-1}\g_1)+
\g_2(\g_2^{-1}\g_1)^{\prime}.
\end{align}
Because of \eqref{3.1}, \eqref{p18} and formula above we get
\begin{align}&      \nn
\g_2^{\prime}\big(I_{m_2}-(\g_2^{-1}\g_1)(\g_2^{-1}\g_1)^*\big)=
\g_2(\g_2^{-1}\g_1)^{\prime}(\g_2^{-1}\g_1)^*, \quad \g_2(0)=I_{m_2}.
\end{align}
Therefore, taking into account that
$I_{m_2}-(\g_2^{-1}\g_1)(\g_2^{-1}\g_1)^*=\g_2^{-1}(\g_2^{-1})^*$
is invertible, we obtain
\begin{align}&      \label{p19}
\g_2^{\prime}=
\g_2(\g_2^{-1}\g_1)^{\prime}(\g_2^{-1}\g_1)^*
\big(I_{m_2}-(\g_2^{-1}\g_1)(\g_2^{-1}\g_1)^*\big)^{-1}, \quad \g_2(0)=I_{m_2}.
\end{align}
\begin{Rk}\label{Schur} If the conditions of Remark \ref{beta'}
hold, then $\g_2$,  $\g$ and Hamiltonian $H=\g^*\g$ are 
consecutively recovered from
$\g_2^{-1}\g_1$ via the differential equation \eqref{p19}. Furthermore,
from the third relation in \eqref{3.8} we see that $(\g_2^{-1}\g_1)(\g_2^{-1}\g_1)^*
<I_{m_2}$. Thus, 
the matrix function $\g_2^{-1}\g_1$ is a continuous analog
of Schur coefficients for the canonical system \eqref{3.14}.
(Compare with Remark 3.1 from \cite{FKS2} on the canonical systems,
which were treated there.)
\end{Rk}
\section{Weyl function: high energy  asymptotics }\label{HE}
\setcounter{equation}{0}
To define Weyl functions, we introduce a class  of nonsingular $m \times m_1$ matrix functions 
$\clp(z)$ with property-$j$, which are an immediate analog of the classical pairs
of parameter matrix functions. Namely, the matrix functions 
$\clp(z)$ are meromorphic in $\BC_+$ and satisfy
(excluding, possibly, a discrete set of points)
the following relations
\begin{align}\label{2.1}&
\clp(z)^*\clp(z) >0, \quad \clp(z)^* j \clp(z) \geq 0 \quad (z\in \BC_+).
\end{align}
It is immediate  from \eqref{1.1} that
\begin{align}&      \label{2.4}
\frac{d}{dx}\big(u(x,z)^*ju(x,z)\big)=\I (z-\ov z)u(x,z)^*u(x,z)<0, \quad z \in \BC_+.
\end{align}
Relations \eqref{2.1} and \eqref{2.4} imply
\begin{align}&      \label{2.7}
\det \Big(\begin{bmatrix}
I_{m_1} & 0
\end{bmatrix}u(x,z)^{-1}\clp(z)\Big)\not= 0,
\end{align}
\begin{Dn} \label{set}
The set $\cln(x,z)$ of M\"obius transformations is the set of values at $x, \,z$ 
of matrix functions
\begin{align}\label{2.2}&
\vp(x,z,\clp)=\begin{bmatrix}
0 &I_{m_2}
\end{bmatrix}u(x,z)^{-1}\clp(z)\Big(\begin{bmatrix}
I_{m_1} & 0
\end{bmatrix}u(x,z)^{-1}\clp(z)\Big)^{-1},
\end{align}
where $\clp(z)$ are nonsingular  matrix functions 
 with property-$j$. 
 \end{Dn}
We can rewrite \eqref{2.2} in an equivalent form, which will be used later on
\begin{align}&      \label{2.9}
\begin{bmatrix}
I_{m_1} \\ \vp(x,z,\clp)
\end{bmatrix}=u(x,z)^{-1}\clp(z)\Big(\begin{bmatrix}
I_{m_1} & 0
\end{bmatrix}u(x,z)^{-1}\clp(z)\Big)^{-1}.
\end{align}

 \begin{Pn}\cite{FKRSp1} \label{PnW1} Let Dirac system \eqref{1.1} on $[0, \, \infty)$
 be given and assume that $v$ is  locally summable.
 Then the sets $\cln(x,z)$
 are well-defined. There is a unique matrix function
 $\vp(z)$ in $\BC_+$ such that
\begin{align}&      \label{2.3}
\vp(z)=\bigcap_{x<\infty}\cln(x,z).
\end{align} 
This function is analytic and non-expansive.
 \end{Pn}
In view of Proposition \ref{PnW1}
we define the Weyl function of Dirac system similar to the canonical system 
case \cite{SaL3}. 
\begin{Dn} \cite{FKRSp1} \label{defWeyl} The Weyl-Titchmarsh (or simply Weyl) function of Dirac system \eqref{1.1} on $[0, \, \infty)$,
where  potential $v$ is locally summable, is the function $\vp$
given by \eqref{2.3}.
\end{Dn}
By Proposition  \ref{PnW1} the Weyl-Titchmarsh function always exists.
Clearly, it is unique.
 
\begin{Cy} \cite{FKRSp1} \label{CyW2} Let the conditions of Proposition \ref{PnW1} hold.
Then the Weyl function is the unique function, which satisfies the inequality
\begin{align}&      \label{2.20}
\int_0^{\infty}
\begin{bmatrix}
I_{m_1} & \vp(z)^*
\end{bmatrix}
u(x,z)^*u(x,z)
\begin{bmatrix}
I_{m_1} \\ \vp(z)
\end{bmatrix}dx< \infty .
\end{align}
\end{Cy}
\begin{Rk}\label{RkDefW} In view of Corollary \ref{CyW2}, inequality
\eqref{2.20} can be used as an equivalent definition of the Weyl function.
Definition of the form \eqref{2.20} is a more classical one and deals with
solutions of \eqref{1.1} which belong to $L^2(0, \, \infty)$. Compare with definitions
of Weyl-Titchmarsh or $M$-functions for discrete and continuous
systems in \cite{CG2, LS, Mar, SaA1, SaA2, SaL3, T1, T2} (see also references therein).
\end{Rk}

Important works  by  F. Gesztesy and B. Simon \cite{GeSi0, GeSi, Si}  
gave rise to a whole series of interesting papers and results
on the high energy asymptotics of Weyl functions and  local Borg-Marchenko-type
uniqueness theorems (see, e.g., \cite{CG1, CGZ, CGZ2, FKS2, LaWo, 
SaA3} and references therein). Here we generalize the high energy
asymptotics result from \cite{SaA3} for the case that Dirac system \eqref{1.1}
has  a rectangular $m_1 \times m_2$ potential $v$, where $m_1$ is not necessarily equal to $m_2$. 
For that we recall first that $S>0$ and $\Phi_1$ is boundedly differentiable
in Lemma \ref{LaNode}.
Therefore, using   Theorem 2.5 from \cite{FKRSp2} we get the statement below.
\begin{Tm}\label{TmIdent}
Let  $\Pi=[\Phi_1 \quad \Phi_2]$ be constructed in Lemma \ref{LaNode}.
Then there is a unique solution 
$S \in B\big(L^2_{m_2}(0,l)\big)$ of the operator
identity \eqref{p9}, this $S$ is strictly positive $($i.e., $S>0)$ and is  defined by the equalities
\begin{align}&
\label{4.5}
\big(Sf\big)(x)=\big(I_{m_2}-\Phi_1(0)\Phi_1(0)^*\big)f(x)-\int_0^ls(x,t)f(t)dt, \\
&\nn
s(x,t):=\int_0^{\min(x,t)}\Phi_1^{\prime}(x-\zeta)\Phi_1^{\prime}(t-\zeta)^*d\zeta
+
\begin{cases}\Phi_1^{\prime}(x-t)\Phi_1(0)^*, \quad x>t;
\\
\Phi_1(0)\Phi_1^{\prime}(t-x)^*, \quad t>x.
\end{cases}
\end{align}
\end{Tm}
Continuous operator kernels of the form above (with a possible jump at
$x=t$) were considered
in Section 2.4 from \cite{gohkol}, where they were called "close to
displacement kernels" (see also references therein).

Now, we will apply the $S$-node scheme, which was used in \cite{SaA1}
for the skew-self-adjoint case,
to derive the high energy asymptotics of  the Weyl function
of Dirac system with a rectangular potential.
\begin{Tm}\label{TmHea} Assume that $\vp \in \cln(l,z)$, where $ \cln(l,z)$ is defined
in Definition \ref{set} and the potential $v$ of the corresponding   Dirac system \eqref{1.1} 
is bounded on $[0, \, l]$. Then $($uniformly with respect to $\Im(z))$ we have
\begin{align} \label{hea}&
\vp(z)=2\I z\int_0^l\E^{2\I xz}\Phi_1(x)dx+O\left(2z \E^{2\I lz}/ \sqrt{\Im(z)}\right), \quad \Im (z) \to \infty.
\end{align}
\end{Tm}
\begin{proof}. To prove the theorem, we consider the matrix function
\begin{align}&\label{4.6}
\clu(z)=\begin{bmatrix}
I_{m_1} & \vp(z)^*
\end{bmatrix}\big(j-
w_A(l,2z)^*jw_A(l,2z)\big)
\begin{bmatrix}
I_{m_1} \\ \vp(z)
\end{bmatrix}.
\end{align}
Because of \eqref{3.16} and \eqref{p17} we have
\begin{align}&\nn
\clu(z)=I_{m_1}-\vp(z)^*\vp(z)- e^{\I l(\ov z-  z)}\begin{bmatrix}
I_{m_1} & \vp(z)^*
\end{bmatrix}
u(l,z)^*ju(l,z)
\begin{bmatrix}
I_{m_1} \\ \vp(z)
\end{bmatrix}.
\end{align}
Taking into account  \eqref{2.9}, we rewrite the formula above as
\begin{align}&\label{4.8}
\clu(z)=I_{m_1}-\vp(z)^*\vp(z)-e^{\I l(\ov z-  z)} \\ & \nn
\times
\Big(\big(\begin{bmatrix}
I_{m_1} & 0
\end{bmatrix}u(x,z)^{-1}\clp(z)\big)^{-1}\Big)^*
\clp(z)^*j\clp(z)
\big(\begin{bmatrix}
I_{m_1} & 0
\end{bmatrix}u(x,z)^{-1}\clp(z)\big)^{-1}.
\end{align}
From \eqref{2.1} and \eqref{4.8} we see that
\begin{align}&\label{4.9}
\clu(z)\leq I_{m_1}.
\end{align}
It easily follows from \eqref{p9} and  \eqref{p13a} (see, e.g., \cite{SaL2, FKS2})
that
\begin{align} \label{4.7}&
w_A(l, z)^*jw_A(l,z)=j+\I(z - \ov z)\Pi^*(I-\ov z A^*)^{-1}S_l^{-1}(I-zA)^{-1}\Pi .
\end{align}
Now, formulas \eqref{4.6}, \eqref{4.9} and \eqref{4.7} imply
\begin{align}&\label{4.10}
 2\I (\ov z-z) \begin{bmatrix}
I_{m_1} & \vp(z)^*
\end{bmatrix}
\Pi^*(I-2\ov z A^*)^{-1}S_l^{-1}(I-2zA)^{-1}\Pi 
\begin{bmatrix}
I_{m_1} \\ \vp(z)
\end{bmatrix}\leq I_{m_1}.
\end{align}
Since $S$ is positive and bounded, inequality \eqref{4.10} yields
\begin{align}&\label{4.11}
\left\| (I-2zA)^{-1}\Pi 
\begin{bmatrix}
I_{m_1} \\ \vp(z)
\end{bmatrix}\right\| \leq C/\sqrt{\Im z} \quad {\mathrm{for}}\,\, 
{\mathrm{some}} \quad C>0.
\end{align}
We easily check directly (see also these formulas in the works on the case
$m_1=m_2$) that
\begin{align} \label{4.12}&
\Phi_2^*(I-2zA)^{-1}f=\int_0^l\E^{2\I(x-l)z}f(x)dx, 
\\  \label{4.13}&
  \Phi_2^*(I-2zA)^{-1}\Phi_2=\frac{\I}{2z}\big(\E^{-2\I lz}-1\big)I_{m_2}.
\end{align}
Because of \eqref{4.11}-\eqref{4.13} (and after applying $-\I\Phi_2^*$ to the operator in the
left-hand side part of \eqref{4.11}), we get
\begin{align} \label{4.14}&
\frac{1}{2z}\big(\E^{-2\I lz}-1\big)\vp(z)=\I \E^{-2\I lz}\int_0^l\E^{2\I xz}\Phi_1(x)dx+
O\left(\frac{1}{\sqrt{\Im(z)}}\right).
\end{align}
Since $\vp$ is non-expansive, we see  from \eqref{4.14} that \eqref{hea} holds.
\end{proof}
Now, consider a potential $v$, which is locally bounded, that is, bounded
on all the finite intervals $[0, \, l]$. The following integral representation
is essential in interpolation and inverse problems.
\begin{Cy}\label{cyHea} Let $\vp$ be the Weyl fuction of Dirac system \eqref{1.1}
on  $[0, \, \infty)$, where the potential $v$    
is locally bounded. Then  we have
\begin{align} \label{repr}&
\vp(z)=2\I z\int_0^{\infty}\E^{2\I xz}\Phi_1(x)dx, \quad \Im (z) >0.
\end{align}
\end{Cy}
\begin{proof}. Since $\vp$ is analytic and non-expansive in $\BC_+$, 
for any $\ve>0$ it admits (see, e.g., \cite[Theorem V]{WP}) a representation
\begin{align} \label{repr0}&
\vp(z)=2\I z\int_0^{\infty}\E^{2\I xz}\Phi (x)dx, \quad \Im (z) >\ve>0,
\end{align}
where $\E^{-2\ve x}\Phi(x) \in L^2_{m_2\times m_1}(0, \, \infty)$. Because of \eqref{hea}
and \eqref{repr0} we get
\begin{align}\nn
\wt Q(z):&=\int_0^{l}\E^{2\I(x-l)z}\big(\Phi_1 (x)-\Phi(x)\big)dx
\\  \label{4.15} &
=\int_l^{\infty}\E^{2\I(x-l)z}\Phi(x)dx+O\big(1/\sqrt{\Im (z)}\big).
\end{align}
From \eqref{4.15} we see that $\wt Q(z)$ is bounded in some semi-plane
$\Im(z)\geq \eta_0>0$. Clearly, $\wt Q(z)$ is bounded also in the semi-plane
$\Im(z)< \eta_0$. Since $\wt Q$ is analytic and bounded in $\BC$
and  tends to zero on some rays, we have
\begin{align} \label{4.16} &
\wt Q(z)=\int_0^{l}\E^{2\I(x-l)z}\big(\Phi_1 (x)-\Phi(x)\big)dx \equiv 0.
\end{align}
It follows from \eqref{4.16} that $\Phi_1(x)\equiv \Phi(x)$ on all the
finite intervals $[0, \, l]$. Hence, \eqref{repr0} implies \eqref{repr}.
\end{proof}
\begin{Rk}\label{RkL}
Since $\Phi_1\equiv \Phi$, we get that $ \Phi_1(x)$ does not depend
on $l$ for $l>x$. Compare this with the proof of Proposition 4.1 in \cite{FKS2},
where the  fact that $E(x,t)$ (and so $\Phi_1$) does not depend on $l$
follows from the uniqueness of the factorizations of operators $S_l^{-1}$.
See also Section 3 in \cite{AGKLS2} on the uniqueness of the accelerant.
Furthermore, since $\Phi_1\equiv \Phi$ the proof of  Corollary \ref{cyHea}
implies also that
$\E^{-\ve x}\Phi_1(x) \in L^2_{m_2\times m_1}(0, \, \infty)$ for any $\ve >0$.
\end{Rk}
\begin{Rk}\label{RkAmp}
From \eqref{repr} we see that $\Phi_1^{\prime}$ is a Dirac system analog
of $A$-amplitude, which was studied in \cite{GeSi, Si}.
On the other hand $\Phi_1^{\prime}$ is closely related to the
so called accelerant, which appears for the case that $m_1=m_2$
in papers by M. Krein (see, e.g., \cite{AGKLS2, Kre1, K}). See also Remark 2.3 in \cite{FKS2},
where $\Phi_1^{\prime}$ is discussed for the case that $m_1=m_2$.
\end{Rk}
\section{Inverse problem and Borg-Marchenko-type uniqueness
theorem} \label{IP}
\setcounter{equation}{0}
Taking into account Plancherel Theorem and Remark \ref{RkL},
we apply inverse Fourier transform to formula \eqref{repr} and get
\begin{align} \label{5.1}&
\Phi_1\Big(\frac{x}{2}\Big)=\frac{1}{\pi}\E^{x\eta}{\mathrm{l.i.m.}}_{a \to \infty}
\int_{-a}^a\E^{-\I x\xi}\frac{\vp(\xi+\I \eta)}{2\I(\xi +\I \eta)}d\xi, \quad \eta >0.
\end{align}
Here l.i.m. stands for the entrywise limit in the norm of  $L^2(0,r)$, 
$\, 0<r \leq \infty$. 
(Note that if we put additionally $\Phi_1(x)=0$ for $x<0$, equality \eqref{5.1}
holds for l.i.m. as the entrywise limit in $L^2(-r,r)$.)
Thus, operators $S$ and $\Pi$ are recovered from $\vp$. 

Since Hamiltonian
$H=\g^*\g$ is recovered from $S$ and $\Pi$ via formula \eqref{p14'}, 
we recover also $\g$. 
Indeed, first we recover  the Schur coefficient (see Remark \ref{Schur} for the motivation
of the term ``Schur coefficient''):                                                                                                                                  
\begin{align} \label{5.2}&
\left(\begin{bmatrix}
0 &I_{m_2}
\end{bmatrix}H\begin{bmatrix}
0 \\ I_{m_2}
\end{bmatrix}\right)^{-1}
\begin{bmatrix}
0 &I_{m_2}
\end{bmatrix}H\begin{bmatrix}
I_{m_1} \\ 0
\end{bmatrix}=(\g_2^*\g_2)^{-1}\g_2^*\g_1=\g_2^{-1}\g_1.
\end{align}
Next, we recover $\g_2$ from $\g_2^{-1}\g_1$ using \eqref{p19}, and finally,
we easily recover $\g$ from $\g_2$ and $\g_2^{-1}\g_1$.

To recover $\b$ from $\g$, we separate $\b$ into two blocks $\b=\begin{bmatrix}
\b_1 & \b_2
\end{bmatrix}$, where $\b_k$ ($k=1,2$) is an $m_1\times m_k$ matrix function.
We put
\begin{align} \label{5.3}&
\wt \b=\begin{bmatrix}
I_{m_1} & \g_1^*(\g_2^*)^{-1}
\end{bmatrix}.
\end{align}
Because of \eqref{3.3} and \eqref{5.3}, we have $\b j \g^*=\wt \b j \g^*=0$, and so
\begin{align} \label{5.4}&
\b(x)= \b_1(x)\wt \b(x).
\end{align}
It follows from \eqref{1.1} and \eqref{3.1} that 
\begin{align} \label{5.4'}&
\b^{\prime}(x)=\I v(x)\g(x),
\end{align}
which implies
\begin{align} \label{5.5}&
\b^{\prime}j\b^*=0, \qquad \b^{\prime}j\g^*=-\I v.
\end{align}
Formula \eqref{5.4} and the first relations in \eqref{3.3} and \eqref{5.5} lead us to
\begin{align} \label{5.6}&
\b^{\prime}j\b^*=\b_1^{\prime}\b_1^{-1}+\b_1(\wt \b^{\prime}j\wt \b^*)\b_1^*=0,
\quad \wt \b j \wt \b^*=\b_1^{-1}(\b_1^*)^{-1}.
\end{align}
By \eqref{1.3} and \eqref{5.6}  $\b_1$ satisfies the first order differential equation
\begin{align} \label{5.7}&
\b_1^{\prime}=-\b_1(\wt \b^{\prime}j\wt \b^*)( \wt \b j \wt \b^*)^{-1}, \quad \b_1(0)=I_{m_1}.
\end{align}
Relations \eqref{5.1}-\eqref{5.5} and \eqref{5.7} give us a procedure to construct the solution
of the inverse problem.
 \begin{Tm} \label{TmIP} Let $\vp$ be the Weyl function of Dirac system \eqref{1.1} on $[0, \, \infty)$,
 where the potential $v$ is  locally bounded.
 Then $v$ can be uniquely recovered from $\vp$ via the formula
\begin{align}&      \label{5.8}
v(x)=\I \b^{\prime}(x)j\g(x)^*.
\end{align} 
Here $\b$ is recovered from $\g$ using \eqref{5.3}, \eqref{5.4} and \eqref{5.7}; $\g$ is recovered
from the Hamiltonian $H$ using \eqref{5.2} and equation \eqref{p19}; the Hamiltonian is given by
\eqref{p14'} and $\Pi$ and $S$ in \eqref{p14'} are expressed via $\Phi_1(x)$ in formulas
\eqref{3.19} and \eqref{4.5}. Finally, $\Phi_1(x)$ is recovered from $\vp$ using \eqref{5.1}.
 \end{Tm}
There is another way to recover $\beta$ and $\gamma$:
\begin{Rk} \label{Rkg} We can recover $\beta$
directly from $\Pi$ and $S$ following the proposition below, 
and then recover $\gamma$ from $\beta$ in the same way that
$\beta$ is recovered from $\gamma$.
\end{Rk}
\begin{Pn} \label{Pnbeta} Let Dirac system \eqref{1.1} on $[0, \, \infty)$
 be given. Assume that $v$ is  locally bounded and $\Pi$ and $S$ are operators
constructed in Lemma \ref{LaNode}. Then the matrix function $\b$, which is defined in \eqref{3.1},
satisfies the equality 
\begin{align}&      \label{5.9}
\b(x)=\begin{bmatrix}I_{m_1} &0 \end{bmatrix}
+\int_0^x\Big(S_x^{-1}\Phi_1^{\prime}\Big)(t)^*\begin{bmatrix}\Phi_1(t) & I_{m_2} \end{bmatrix}dt.
\end{align} 
 \end{Pn}
\begin{proof}. First, we fix an arbitrary $l$ and rewrite \eqref{p10} (for $x<l$) in the form 
 \begin{align}&      \label{5.10}
\g(x)=\big(E \begin{bmatrix}\Phi_1 & I_{m_2} \end{bmatrix}\big)(x).
\end{align} 
It follows that
 \begin{align}&      \label{5.11}
\I \big(E AE^{-1}E\Phi_1^{\prime}\big)(x)=\g_1(x)-\g_2(x)\Phi_1(+0).
\end{align}
Recalling that $\g_1(0)=0$ and taking into account \eqref{3.19} and Remark \ref{kern},
we have $\Phi_1(+0)=0$. Therefore we rewrite \eqref{5.11} as
 \begin{align}&      \label{5.12}
\I \big(E AE^{-1}E\Phi_1^{\prime}\big)(x)=\g_1(x).
\end{align}
Next, we substitute $K=EAE^{-1}$ from \eqref{3.6} into \eqref{5.12} and (using  \eqref{3.4} and \eqref{3.8}) we  get
 \begin{align}&      \label{5.13}
\g_1(x)=-\g(x)j\int_0^x \g(t)^*\Big(E\Phi_1^{\prime}\Big)(t)dt.
\end{align}
Formulas \eqref{5.10} and \eqref{5.13} imply
 \begin{align}&      \label{5.14}
\g_1(x)=-\g(x)j\int_0^x \big(E \begin{bmatrix}\Phi_1 & I_{m_2} \end{bmatrix}\big)(t)^*\Big(E\Phi_1^{\prime}\Big)(t)dt.
\end{align}
Because of  \eqref{p14!}, \eqref{p14!!} and \eqref{5.14} we see that
\begin{align}     \label{5.15}
\g(x)j\t(x)^*\equiv 0, \quad \t(x):&=\begin{bmatrix}I_{m_1} &0 \end{bmatrix}+\int_0^x 
\Big(E\Phi_1^{\prime}\Big)(t)^*
\big(E \begin{bmatrix}\Phi_1 & I_{m_2} \end{bmatrix}\big)(t)dt \\
&\nn =
\begin{bmatrix}I_{m_1} &0 \end{bmatrix}
+\int_0^x\Big(S_x^{-1}\Phi_1^{\prime}\Big)(t)^*\begin{bmatrix}\Phi_1(t) & I_{m_2} 
\end{bmatrix}dt,
\end{align}
where $S_x:=P_xSP_x$. We shall show that $\t=\b$.

In view of \eqref{5.10} and the second equality in \eqref{5.15},
we have 
\begin{align}     \label{5.15'}&
\t^{\prime}(x)=\Big(E\Phi_1^{\prime}\Big)(x)^*\g(x).
\end{align}
Therefore, \eqref{3.3}
leads us to
\begin{align}     \label{5.16}&
\b(x)j\t^{\prime}(x)^* \equiv 0.
\end{align}
Furthermore, compare \eqref{3.3} with the first  equality in \eqref{5.15}
to see that
\begin{align}     \label{5.17}&
\t(x)=\vk(x)\b(x),
\end{align}
where $\vk(x)$ is an $m_1 \times m_1$ matrix function, which is boundedly differentiable
on $[0, \, l]$. Now, equalities \eqref{5.16} and \eqref{5.17}
and the first relations in \eqref{3.3} and \eqref{5.5} yield that $\vk^{\prime}\equiv 0$
(i.e., $\vk$ is a constant). It follows from \eqref{3.1}, \eqref{5.15} and \eqref{5.17}
that $\vk(0)=I_{m_1}$, and so $\vk \equiv I_{m_1}$, that is, $\t \equiv \b$.
Thus, \eqref{5.9} is immediate from \eqref{5.15}.
\end{proof}
\begin{Rk} Because of \eqref{5.4'}, \eqref{5.15'} and equality $\b=
$, we see that
the potential $v$ can be recovered via the formula
\begin{align}     \label{5.18}&
v(x)=\Big(\I E\Phi_1^{\prime}\Big)(x)^*.
\end{align}
Formally applied, formulas \eqref{p14!!}, \eqref{3.10} and \eqref{5.18}
yield
\begin{align}     \label{5.19}&
v(x)=\Big(\I S_x^{-1}\Phi_1^{\prime}\Big)(x)^*,
\end{align}
though one needs a proper "pointwise" definition of  matrix functions
$S_x^{-1}\Phi_1^{\prime}$ for \eqref{5.19} to hold.
\end{Rk}
The last statement in this section is a Borg-Marchenko-type uniqueness
theorem, which follows from Theorems \ref{TmHea} and  \ref{TmIP}.

\begin{Tm}\label{BM} Let $\vp$ and $\wh \vp$ be the Weyl functions
of two Dirac systems on $[0, \, \infty)$ with the locally bounded potentials,
which are denoted by $v$ and $\wh v$, respectively. Suppose that on some ray
$\Re z=c \Im z$ $\,(c \in \BR, \, \Im z>0)$ the equality
\begin{align}     \label{5.20}&
\|\vp(z)-\wh \vp(z)\|=O(\E^{2\I rz}) \quad (\Im z \to \infty) \quad {\mathrm{for}}
\,\, {\mathrm{all}} \quad 0<r<l
\end{align}
holds. Then we have
\begin{align}     \label{5.21}&
v(x)=\wh v(x), \qquad 0<x<l.
\end{align}
\end{Tm}
\begin{proof}.  Since Weyl functions are non-expansive, we get
\begin{align}     \label{5.22}&
\|\E^{-2\I rz}\big(\vp(z)-\wh \vp(z)\big)\|\leq c_1 \E^{2r|z|}, \quad \Im z \geq c_2>0
\end{align}
for some $c_1$ and $c_2$, and
the matrix function $\E^{-2\I rz}\big(\vp(z)-\wh \vp(z)\big)$ is bounded
on  the line $\Im z =c_2$. Furthermore, formula \eqref{5.20} implies that
$\E^{-2\I rz}\big(\vp(z)-\wh \vp(z)\big)$ is bounded on the ray $\Re z=c \Im z$.
Therefore, applying the Phragmen-Lindel\"of theorem in the angles
generated by the line $\Im z =c_2$ and the ray $\Re z=c \Im z$ ($\Im z \geq c_2$),
we see that
\begin{align}     \label{5.23}&
\|\E^{-2\I rz}\big(\vp(z)-\wh \vp(z)\big)\|\leq c_3, \quad \Im z \geq c_2>0.
\end{align}
Functions associated with $\wh \vp$ will be written with a hat (e.g.,
$\wh \Phi_1$).  Because of formula \eqref{hea}, its analog for $\wh \vp$,
$\wh \Phi_1$ and the inequality \eqref{5.23}, we have
\begin{align}     \label{5.24}&
\| \int_0^r \E^{2\I(x- r)z}\big(\Phi_1(x)-\wh \Phi_1(x)\big)dx   \|
\leq c_4, \quad \Im z \geq c_2>0.
\end{align}
Clearly, the left-hand side of \eqref{5.24} is bounded in the semi-plane
$\Im z<c_2$ and tends to zero on some rays. Thus, we derive
\begin{align}     \label{5.25}&
\int_0^r \E^{2\I(x- r)z}\big(\Phi_1(x)-\wh \Phi_1(x)\big)dx   \equiv 0,
\quad {\mathrm{i.e.}}, \quad  \Phi_1(x)\equiv \wh \Phi_1(x) \quad
(0<x<r).
\end{align}
Since \eqref{5.25} holds for all $r<l$, we obtain $ \Phi_1(x)\equiv \wh \Phi_1(x)$ for $0<x<l$. In view of Theorem \ref{TmIP}, the last identity implies \eqref{5.21}.
\end{proof}

\section{Weyl function and positivity of $S$} \label{Sc}
\setcounter{equation}{0}
In this section we discuss some sufficient conditions for a non-expansive
matrix function $\vp$, which is analytic in $\BC_+$, to be a Weyl function of 
Dirac system on the semi-axis. For the case that $v$ is a scalar and Weyl
functions $\vp$ are so called Nevanlinna functions (i.e., $\Im \vp \geq 0$),
sufficient condition for $\vp$ to be a Weyl function can be given in terms
of spectral function \cite{LS, Mar},  which is connected with $\vp$ via Herglotz
representation. Furthermore, a positive operator $S$ is also recovered from the spectral function (see \cite[Chapters 4,10]{SaL3} and \cite{SaA3}). 
The invertibility of the convolution operators, which is required in
 \cite{AGKLS1, AGKLS2, Kre1}, provides their positivity too, and the spectral
 problem is treated in this way.
Here, we again formulate conditions on the $m_2 \times m_1$ non-expansive
matrix functions in terms of $S$.
To derive our sufficient conditions, we apply the procedure to recover
Dirac system from its Weyl function (see Section \ref{IP}).
Recall that $\Phi_1$ in Section \ref{IP} is the Fourier transform of $\vp$, that is, it  is given by \eqref{5.1}.

First, consider a useful procedure to recover $\g$ from $\b$ as mentioned in Remark \ref{Rkg}. 
\begin{Pn} \label{PnV} Let a given $m_1 \times m$ matrix function $\b(x)$ $($$0 \leq x \leq l$$)$
be boundedly differentiable and satisfy relations
\begin{align}&      \label{6.11}
\b(0)=\begin{bmatrix}
I_{m_1} & 0
\end{bmatrix}, \qquad \b^{\prime}j\b^* \equiv 0.
\end{align} 
Then there is a unique  $m_2 \times m$ matrix function $\g$,
which is boundedly differentiable  and satisfies relations
\begin{align}&      \label{6.11'}
\g(0)=\begin{bmatrix}
0 & I_{m_2}
\end{bmatrix}, \quad \g^{\prime}j\g^* \equiv 0, \quad  \g j\b^*\equiv 0.
\end{align} 
Moreover, this $\g$ is given by the formula
\begin{align}&      \label{6.11!}
\g=\g_2\wt \g, \quad \wt \g=\begin{bmatrix}
\wt \g_1 & I_{m_2}
\end{bmatrix}, \quad \wt \g_1=\b_2^*(\b_1^*)^{-1},
\end{align} 
where $\det \b_1(x)\not=0$, $\det\big(I_{m_2}- \wt \g_1(x)  \wt \g_1(x)^*\big)\not=0$ and $\g_2$ can be recovered via  the linear differential system and initial
condition below:
\begin{align}&      \label{6.11!!}
\g_2^{\prime}=\g_2\wt \g_1^{\prime} \wt \g_1^*\big(I_{m_2}- \wt \g_1  \wt \g_1^*\big)^{-1}, \quad \g_2(0)=I_{m_2}.
\end{align}
 \end{Pn}
 \begin{proof}.  Because of \eqref{6.11} we have $\b j \b^* \equiv I_{m_1}$ (and so $\det \b_1 \not=0$).  
On the other hand, formula
\eqref{6.11!} implies $\wt \g j \b^*=0$. Therefore, we get
  \begin{align}&      \label{6.12}
\g j \b^*=0, \qquad  \wt \g j \wt \g^*<0.
\end{align} 
In particular, we see that $\det\big(I_{m_2}- \wt \g_1  \wt \g_1^*\big)\not=0$. 
According to \eqref{6.11!!} we have $\det \g_2 \not= 0$.
Furthermore,
any $\g$ satisfying \eqref{6.11'} admits representation $\g=\g_2\wt \g$ with some
boundedly differentiable $\g_2$ such that $\g_2(0)=I_{m_2}$
and  $\det \g_2(x) \not=0$. Thus, it remains to rewrite
$\g^{\prime} j \g^*$ in the equivalent form
 \begin{align}&      \label{6.12'}
\big(\g_2^{\prime}\wt \g+\begin{bmatrix}
\g_2\wt \g_1^{\prime} & 0
\end{bmatrix}\big)j \wt \g^*=0,
\end{align} 
which, in its turn, is equivalent to the first equality in \eqref{6.11!!}. Clearly \eqref{6.11!!} uniquely defines
$\g_2$.
 \end{proof}
Now, we formulate the main statement in this section.
\begin{Tm}\label{TmS}
Let  an $m_2 \times m_1$ matrix function $\vp(z)$ be analytic and non-expansive in $\BC_+$.
Furthermore, let matrix function $\Phi_1(x)$ and  operators $S_l$, which are given by  \eqref{5.1}
and \eqref{4.5}, respectively, be such that $\Phi_1$
is boundedly  differentiable  on each finite interval $[0,\, l]$ and satisfies
equality $\Phi_1(x)=0$ for $x \leq 0$,
whereas operators  $S_l$ 
are boundedly invertible for all $\, 0<l < \infty$. 

Then $\vp$ is the Weyl function
of some Dirac system on $[0, \, \infty)$. 
The operators $S_l^{-1}$ admit unique factorizations
\begin{align}&      \label{6.13}
S_l^{-1}=E_{\Phi,l}^*E_{\Phi,l}, \quad E_{\Phi,l}=I+\int_0^x E_{\Phi}(x,t)\, \cdot \, dt \in
B\big(L_2^{m_2}(0,l)\big),
\end{align} 
where $E_{\Phi}(x,t)$ is continuous with respect to $x, \, t$ and does not depend
on $l$, and
the potential of the Dirac system is constructed via formula
 \begin{align}&      \label{6.14}
v(x)=\Big(\I E_{\Phi,l} \Phi_1^{\prime}\Big)(x)^*, \quad 0<x<l.
\end{align} 
\end{Tm}  

Note that formula \eqref{6.14} in Theorem \ref{TmS}
is similar to \eqref{5.18}. 
To prove this theorem we need  an auxiliary Proposition 2.1 from \cite{FKRSp2}
on operator $S$:
\begin{Pn}\cite{FKRSp2}\label{OpI} Let $\Phi_1(x)$ be an $m_2 \times m_1$ matrix
function, which is boundedly differentiable on the interval $[0, \, l]$.
Then the operator $S$, which is given by \eqref{4.5}, satisfies the
operator identity \eqref{p9}, where $\Pi:= \begin{bmatrix}
\Phi_1 & \Phi_2
\end{bmatrix}$ is expressed via formulas
\begin{align}&      \label{6.2}
\Phi_k \in B\big(\BC^{m_k}, \, L^2_{m_2}(0, \, l)\big),
\quad
\big(\Phi_1 f\big)(x)=\Phi_1(x)f,  \quad  \Phi_2 f=I_{m_2}f.
\end{align}
\end{Pn}
Factorization results from \cite[pp. 185-186]{GoKr} yield  the following lemma (see  \cite[Theorem 4.2]{FKRSp2}).
\begin{La} \label{LaFact}
Let matrix function $\Phi_1(x)$ and  operators $S_l$, which are expressed via $\Phi_1$
in \eqref{4.5},  be such that $\Phi_1$
is boundedly   differentiable on each finite interval $[0,\, l]$ and satisfies
equality $\Phi_1(0)=0$,
whereas operators  $S_l$ 
are boundedly invertible for all $\, 0<l < \infty$. Then the operators $S_l^{-1}$ admit factorizations \eqref{6.13},
where $E_{\Phi}(x,t)$ is continuous with respect to $x, \, t$ and does not depend
on $l$. Furthermore, all the factorizations \eqref{6.13} with continuous $E_{\Phi}(x,t)$
are unique.
\end{La}
Now, we consider our procedure to solve inverse problem.
\begin{La} \label{bega} Let the conditions of Lemma \ref{LaFact} hold.
Then the matrix functions
\begin{align}&      \label{6.16}
\b_{\Phi}(x):=\begin{bmatrix}I_{m_1} &0 \end{bmatrix}
+\int_0^x\Big(S_x^{-1}\Phi_1^{\prime}\Big)(t)^*\begin{bmatrix}\Phi_1(t) & I_{m_2} \end{bmatrix}dt,
\\ &      \label{6.17} \g_{\Phi}(x):=\begin{bmatrix}\Phi_1(x) & I_{m_2} \end{bmatrix}
+\int_0^xE_{\Phi}(x,t)\begin{bmatrix}\Phi_1(t) & I_{m_2} \end{bmatrix}dt
\end{align} 
are boundedly differentiable and satisfy conditions
\begin{align}&      \label{6.18}
\b_{\Phi}(0):=\b_{\Phi}(+0)=\begin{bmatrix}I_{m_1} &0 \end{bmatrix}, \quad
\b_{\Phi}^{\prime}j\b_{\Phi}^*\equiv 0; 
\\ &      \label{6.19} 
\g_{\Phi}(0)=\begin{bmatrix}0& I_{m_2} \end{bmatrix}, \quad
\g_{\Phi}^{\prime}j\g_{\Phi}^*\equiv 0;  
\quad \g_{\Phi}j\b_{\Phi}^*\equiv 0.
\end{align} 
\end{La}
\begin{proof}.  
Step 1.
The first equalities in \eqref{6.18} and  \eqref{6.19}
are immediate from  \eqref{6.16} and  \eqref{6.17}, respectively.
Furthermore, \eqref{6.17} is equivalent to the equalities
\begin{align}&      \label{6.20}
 \g_{\Phi}(x)=\Big(E_{\Phi,l}\begin{bmatrix}\Phi_1 & I_{m_2} \end{bmatrix}\Big)(x),
 \quad 0 \leq x \leq l \quad ({\mathrm{for}} \,\, {\mathrm{all}} \,\, l<\infty).
\end{align} 
Next, fix any $0<l<\infty$ and recall that according to Proposition \ref{OpI}
the operator identity
\begin{align}&      \label{6.21}
 AS_l-S_lA^*=\I \Pi j \Pi^*
 \end{align} 
 holds, where $\Pi$ is expressed via \eqref{6.2}. Hence, taking into account
 \eqref{6.13} and \eqref{6.20}, and turning around the proof of \eqref{6.21} in Lemma \ref{LaNode}, we get
\begin{align}&\nn
E_{\Phi}AE_{\Phi}^{-1}-\big(E_{\Phi}^{-1}\big)^*A^*E_{\Phi}^*=\I\g_{\Phi}(x)j\int_0^l\g_{\Phi}(t)^*\,\cdot \,dt,
\,\,{\mathrm{i.e.}},
 \\ &      \label{6.22}
 E_{\Phi}AE_{\Phi}^{-1}=\I\g_{\Phi}(x)j\int_0^x\g_{\Phi}(t)^*\,\cdot \,dt
 \quad ({\mathrm{for}} \,\, E_{\Phi}=E_{\Phi,l}).
 \end{align} 
Introducing the resolvent kernel $\G_{\Phi}$ by $E_{\Phi}^{-1}=I+\int_0^x\G_{\Phi}(x,t)\cdot dt$, we rewrite \eqref{6.22} in the form of an equality for kernels:
\begin{align}      \nn
I_{m_2}+\int_t^x\big(E_{\Phi}(x,r)+\G_{\Phi}(r,t)\big)dr+
\int_t^x\int_{\xi}^x & E_{\Phi}(x,r)dr\G_{\Phi}(\xi,t)d\xi
\\ &  \label{6.23}
=-\g_{\Phi}(x)j\g_{\Phi}(t)^*.
 \end{align} 
 In particular, formula \eqref{6.23} for the case that $x=t$ implies
 \begin{align}      &  \label{6.24}
 \g_{\Phi}(x)j\g_{\Phi}(x)^* \equiv - I_{m_2}.
 \end{align} 
 In a way quite similar to the first part of the proof of Proposition \ref{Pnbeta},
 we use equalities \eqref{6.20}, \eqref{6.22} and $\Phi_1(0)=0$ to derive
  \begin{align}      &  \label{6.25}
 \g_{\Phi}(x)j\b_{\Phi}(x)^* \equiv 0
 \end{align}
(compare with \eqref{5.15}). Because of \eqref{6.13}, \eqref{6.16} and \eqref{6.20}
we have
\begin{align}&      \label{6.26}
\b_{\Phi}^{\prime}(x)=
\Big(E_{\Phi}\Phi_1^{\prime}\Big)(x)^*\g_{\Phi}(x).
 \end{align}
 In view of  \eqref{6.25} and \eqref{6.26}, $\b_{\Phi}$ is boundedly
 differentiable and the last relation in \eqref{6.18} holds.
 
 Step 2. It remains to show that $\g_{\Phi}$ is boundedly differentiable
 and the identity $\g_{\Phi}^{\prime}j\g_{\Phi}^*\equiv 0$ holds.
 For that purpose note that $\b_{\Phi}$ satisfies conditions of Proposition
 \ref{PnV}, and so there is a boundedly differentiable
 matrix function $\wh \g$ such that
 \begin{align}&      \label{6.27}
\wh \g(0)=\begin{bmatrix}0& I_{m_2} \end{bmatrix}, \quad
\wh \g^{\prime}j \wh \g^*\equiv 0;  
\quad \wh \g j\b_{\Phi}^*\equiv 0.
\end{align} 
Formulas \eqref{6.18} and \eqref{6.27}  yield
  \begin{align}&      \label{6.28}
 \b_{\Phi}j\b_{\Phi}^*\equiv I_{m_1}, \quad 
\wh \g j \wh \g^*\equiv - I_{m_2},
 \end{align} 
 respectively. 
 That is, the rows of $\b_{\Phi}$  (the rows of $\wh \g$)
are linearly independent.  Hence, the
last relations in \eqref{6.27} and \eqref{6.28} and formulas \eqref{6.24}
and \eqref{6.25} imply
that there is a unitary matrix function $\om$ such that
 \begin{align}&      \label{6.29}
 \g_{\Phi}(x)=\om(x) \wh \g(x), \quad \om(x)^*= \om(x)^{-1}.
 \end{align} 
 Moreover, formulas \eqref{6.18},  \eqref{6.27} and \eqref{6.28}
 lead us to the relations
  \begin{align}&      \label{6.30}
\wh u^{\prime}j  \wh u^* j=ij
\begin{bmatrix}0& \wh v \\ \wh v^* &0\end{bmatrix}, \quad
\wh u(x):= \begin{bmatrix}\b_{\Phi}(x) \\ \wh \g(x) \end{bmatrix}, \quad
\wh v:=\I \b_{\Phi}^{\prime}j\wh \g^*;
\\ &      \label{6.31}
\wh u j \wh u^* j \equiv I_m, \quad \wh u(0)=I_m.
\end{align} 
According to \eqref{6.30} and  \eqref{6.31} the matrix function
$\wh u(x)$ is the normalized by \eqref{1.3} solution of Dirac system, 
where $\wh v$ is the bounded potential and the spectral parameter $z$
equals zero. In other words, $\b_{\Phi}$ and $\wh \g$ correspond to
$\wh v$ via equalities \eqref{3.1}, and we can apply the results
of Section \ref{Snode}. Therefore, by Proposition \ref{PnSimN} and
Remark \ref{beta'} there is an operator $E$ of the form \eqref{3.10},
such that
 \begin{align}&      \label{6.32}
EA=i\wh \g(x)j\int_0^x\wh \g(t)^* \, \cdot \, dt\, E, \quad \wh \g_2=EI_{m_2}.
 \end{align} 
 It follows from \eqref{6.29} and \eqref{6.32} that
 \begin{align}&      \label{6.32'}
\om EA=i \g_{\Phi}(x) j\int_0^x \g_{\Phi}(t)^* \, \cdot \, dt\, \om E, \quad  
\g_{2,\Phi}=\om EI_{m_2},
 \end{align}  
 where $\om$ denotes also the operator of multiplication by the matrix function $\om$.
On the other hand, formulas \eqref{6.20} and \eqref{6.22} lead us to
 \begin{align}&      \label{6.33}
E_{\Phi}A=i \g_{\Phi}(x)j\int_0^x \g_{\Phi}(t)^* \, \cdot \, dt \, E_{\Phi}, \quad
 \g_{2,\Phi}=E_{\Phi}I_{m_2}.
 \end{align} 
It is easy to see that c.l.s.$\bigcup_{i=0}^{\infty}\im \big(A^i I_{m_2}\big)=
L^2_{m_2}(0,l)$. Hence, equalities  \eqref{6.32'} and \eqref{6.33}
imply $E_{\Phi}=\om E$. Taking into account that by \eqref{3.10}
and \eqref{6.13} expressions for both operators $E$ and $E_{\Phi}$  
include a term $I$, we see that
 \begin{align}&      \label{6.34}
\om(x) \equiv I_{m_2}, \quad E_{\Phi}=E.
 \end{align} 
 Furthermore, because of \eqref{6.29} and \eqref{6.34} we have,
 $\g_{\Phi}=\wh \g$, and so $\g_{\Phi}$
 is boundedly differentiable and satisfies \eqref{6.19}.
\end{proof}
\begin{proof} of Theorem \ref{TmS}. By the assumptions of Theorem \ref{TmS}
the conditions of Lemmas \ref{LaFact}, \ref{bega} are fulfilled. Therefore, 
the theorem's statements about operators $S_l$ and $E_{\Phi,l}$ are true.
Furthermore, because of \eqref{6.24}, \eqref{6.26}, \eqref{6.30} and
equality $\g_{\Phi}=\wh \g$ we have $v=\wh v$ for $v$ given by \eqref{6.14}
and $\wh v$ from Lemma \ref{bega}. It follows that $v$ is bounded on $[0, \, l]$,
$\g_{\Phi}=\wh \g=\g$
and operator $E$, which is recovered from $v$ in Section \ref{Snode},
satisfies \eqref{6.34}. Thus, we derived:
 \begin{align}&      \label{6.35}
v=\wh v, \quad \g_{\Phi}=\g, \quad E_{\Phi}=E.
 \end{align} 
Using \eqref{3.19}, \eqref{6.20} and \eqref{6.35} we see that
$\Phi_1$, which is recovered from $v$ and $E$ in Section \ref{Snode}, 
coincides with $\Phi_1$ in the statement of theorem. 
In view of Proposition \ref{PnW1} and Corollary \ref{cyHea}
there is a unique Weyl function $\vp_W$ of the  Dirac system
with the potential $v$ and this Weyl function is given by \eqref{repr}.

It remains to show that $\vp_W$ equals the function $\vp$, which generates $\Phi_1$
via \eqref{5.1}. Recalling that \eqref{5.1} holds also for $\vp_W$,
we see that
\begin{align} \label{6.36}&
{\mathrm{l.i.m.}}_{a \to \infty}
\int_{-a}^a\frac{\E^{-\I x\xi}}{(\xi +\I \eta)}\big(\vp(\xi+\I \eta)
-\vp_W(\xi+\I \eta)\big)
d\xi \equiv 0, \quad \eta >0,
\end{align}
where l.i.m. stands for the entrywise limit in the norm of  $L^2(-r,r)$
($\, 0<r \leq \infty$). Therefore, we get $\vp_W=\vp$, that is, 
$\vp$ is the Weyl function of the Dirac system, where  the potential is given by
\eqref{6.14}.
\end{proof}

{\bf Acknowledgement.}
The work of I.Ya. Roitberg was supported by the 
German Research Foundation (DFG) under grant no. KI 760/3-1 and
the work of A.L. Sakhnovich was supported by the Austrian Science Fund (FWF) under Grant  no. Y330.



\begin{flushright} \it
B. Fritzsche,  \\
Fakult\"at f\"ur Mathematik und Informatik, \\
Mathematisches Institut, Universit\"at Leipzig, \\
Johannisgasse 26,  D-04103 Leipzig, Germany,\\
e-mail: {\tt  Bernd.Fritzsche@math.uni-leipzig.de } \\   $ $ \\

B. Kirstein, \\
Fakult\"at f\"ur Mathematik und
Informatik, \\
Mathematisches Institut, Universit\"at Leipzig,
\\ Johannisgasse 26,  D-04103 Leipzig, Germany, \\
e-mai: {\tt Bernd.Kirstein@math.uni-leipzig.de } \\  $ $ \\

I. Roitberg, \\
Fakult\"at f\"ur Mathematik und
Informatik, \\
Mathematisches Institut, Universit\"at Leipzig, \\
Johannisgasse 26,  D-04103 Leipzig, Germany, \\
e-mail: {\tt i$_-$roitberg@yahoo.com } \\  $ $ \\

A.L. Sakhnovich, \\  Fakult\"at f\"ur Mathematik,
Universit\"at Wien,
\\
Nordbergstrasse 15, A-1090 Wien, Austria \\
e-mail: {\tt al$_-$sakhnov@yahoo.com }
\end{flushright}

\end{document}